\documentclass[12pt, reqno]{amsart}
\usepackage{amsmath, amsthm, amscd, amsfonts, amssymb, graphicx, color}
\usepackage[bookmarksnumbered, colorlinks, plainpages]{hyperref}
\hypersetup{colorlinks=true,linkcolor=red, anchorcolor=green, citecolor=cyan, urlcolor=red, filecolor=magenta, pdftoolbar=true}

\def \R {\mathbb{R} }
\newtheorem{theorem}{Theorem}[section]
\newtheorem{lemma}[theorem]{Lemma}
\newtheorem{proposition}[theorem]{Proposition}
\newtheorem{corollary}[theorem]{Corollary}
\theoremstyle{definition}
\newtheorem{definition}[theorem]{Definition}

\newtheorem{remark}[theorem]{Remark}
\numberwithin{equation}{section}

\begin{document}

\setcounter{page}{1}

\title[Existence of solutions]{Existence of solutions for a nonlocal type problem in fractional Orlicz Sobolev spaces}

\author[E. AZROUL, A. BENKIRANE, M. SRATI]{E. AZROUL,  A. BENKIRANE, \MakeLowercase{and} M. SRATI$^{*}$}
\address{
Sidi Mohamed Ben Abdellah
 University,
 Faculty of Sciences Dhar El Mahraz, Laboratory of Mathematical Analysis and Applications, Fez, Morocco.}

\email{\textcolor[rgb]{0.00,0.00,0.84}{ elhoussine.azroul@gmail.com}}
\email{\textcolor[rgb]{0.00,0.00,0.84}{abd.benkirane@gmail.com}}
\email{\textcolor[rgb]{0.00,0.00,0.84}{srati93@gmail.com}}


\let\thefootnote\relax\footnote{Copyright 2016 by the Tusi Mathematical Research Group.}



\date{Received: xxxxxx; Accepted: zzzzzz.
\newline \indent $^{*}$Corresponding author
\newline $^\diamond$ Advance publication -- final volume, issue, and page numbers to be assigned.}

\begin{abstract}
In this paper, we investigate the existence of weak solution for a fractional type problems driven by a nonlocal  operator of elliptic type in a fractional Orlicz-Sobolev space, with homogeneous Dirichlet boundary conditions.
 We first  extend the fractional Sobolev spaces  $W^{s,p}$ to include the general  case $W^sL_A$, where $A$ is an N-function and $s\in (0,1)$. We are concerned with some qualitative properties of the space $W^sL_A$ (completeness, reflexivity and separability). Moreover,  we prove a continuous and compact embedding theorem of these spaces into Lebesgue spaces.\\
\end{abstract} \maketitle

\section{Introduction}
In this paper,  we establish the existence of a weak solutions for  the following Dirichlet type equation  
            \begingroup\makeatletter\def\f@size{10}\check@mathfonts
             $$
           \label{P_a}  (P_a) \hspace*{0.5cm} \left\{ 
               \begin{array}{clclc}
            
               (-\Delta)^s_a u & = & f(x,u)   & \text{ in }& \Omega, \\\\
                u & = & 0  & \text{ in } & \R^N\setminus \Omega,
                \label{eq1}
               \end{array}
               \right. 
            $$ 
            \endgroup
             where $\Omega$ is an open bounded subset in $\R^N$ with Lipschitz boundary $\partial \Omega$, $0<s<1$,   
           $f: \Omega\times \R \longrightarrow \R$ is a Carath\'eodory function and $(-\Delta)^s_a$ is the fractional $a$-Laplacian operator  defined as                                   
                    \begingroup\makeatletter\def\f@size{10}\check@mathfonts
                       $$
                        \begin{aligned}
                        (-\Delta)^s_au(x)&=2p.v \int_{\R^N} A'\left( \dfrac{|u(x)-u(y)|}{|x-y|^s }\right)\dfrac{u(x)-u(y)}{|u(x)-u(y)|}\dfrac{dy}{|x-y|^{N+s}}\\
                        &=2\lim\limits_{\varepsilon\searrow 0} \int_{\R^N\setminus B_\varepsilon(x)} a\left( \dfrac{|u(x)-u(y)|}{|x-y|^s }\right)\dfrac{u(x)-u(y)}{|u(x)-u(y)|} \dfrac{dy}{|x-y|^{N+s} }
                        \end{aligned}
                         $$\endgroup
     with $p.v$. stands for in principal value,  $a=A'$ and  $A$ is an N-function.

 The study of nonlinear elliptic equations involving quasilinear homogeneous type operators is based on the theory of Sobolev spaces and fractional Sobolov spaces $W^{s,p}(\Omega)$  in order to find weak solutions. In certain equations, precisely in the case of nonhomogeneous differential operators, when trying to relax some conditions on these operators (as growth conditions), the problem can not be formulated with classical Lebesgue and Sobolev spaces. Hence, the adequate functional spaces is the so-called Orlicz spaces. These spaces consists of functions that have weak derivatives and satisfy certain integrability conditions. Many properties of Orlicz-Sobolev spaces  come in \cite{1,21,28,34}. For this, many researchers have studied the existence of solutions for the eigenvalue problems involving nonhomogeneous operators in the divergence form through Orlicz-Sobolev spaces by using variational methods and critical point theory, monotone operator methods, fixed point theory and degree theory (see, for instance, \cite{af,3,16,bo1,bo2}).
  
    The problem \hyperref[P_a]{$(P_a)$} involves the fractional $a$-Laplacian  operator, the most appropriate functional framework for dealing with this problem is the fractional Orlicz Sobolev space which introduced by Salort et al \cite{9}, namely a fractional Sobolev space constructed from an Orlicz space at the place of $L^p(\Omega)$. As we know, the Orlicz spaces represent a generalization of classical Lebesgue spaces in which the role usually played by the convex function $t^p$ is assumed by a more general convex function $A(t)$; they have been extensively studied in the monograph of Krasnosel\'skii and Rutickii \cite{28} as well as in Luxemburg's doctoral thesis \cite{Le}. If the role played by $L^p(\Omega)$ in the definition of fractional Sobolev spaces $W^{s,p}(\Omega)$ is assigned to an Orlicz $L_A(\Omega)$ space, the resulting space $W^sL_A(\Omega)$ is exactly a fractional Orlicz-Sobolev space. Many properties of fractional Sobolev spaces have been extended to fractional Orlicz - Sobolev spaces (see section 3).
    
    In applied PDE, fractional spaces and the corresponding nonlocal equations, are now experiencing impressive applications in different subjects, such as, among others, the thin obstacle problem \cite{31}, finance \cite{17}, phase transitions  \cite{2,10}, stratified
    materials \cite{15}, crystal dislocation \cite{7}, soft thin films \cite{29}, semipermeable membranes and flame
    propagation \cite{11}, conservation laws \cite{8}, ultra-relativistic limits of quantum mechanics \cite{26}, quasi-geostrophic flows \cite{13}, multiple scattering \cite{22}, minimal surfaces \cite{12} , materials science \cite{5}, water
    waves \cite{37},
    gradient potential theory \cite{33} and singular set of minima of variational
    functionals \cite{32}. See also \cite{35} for further
    motivation.

  When $A(t)=\dfrac{t^p}{p}$, the problem  \hyperref[P_a]{$(P_a)$} reduces to the fractional $p$-Laplacian problem
  $$
 \label{P}  (P_p)\hspace*{0.3cm} \left\{ 
     \begin{array}{clclc}
   (-\Delta)^s_pu & = & f(x,u)   & \text{ in }& \Omega \\\\
     \hspace*{1.5cm} u & = & 0 \hspace*{0.2cm} \hspace*{0.2cm} & \text{ in } & \R^N\smallsetminus \Omega,
     \end{array}
     \right. 
  $$
  where $(-\Delta)^s_p$ is the fractional $p$-Laplacian operator which, up to normalization, may defined as
  $$(-\Delta)^s_p u(x)=2\lim_{\varepsilon \searrow 0}\int_{\R^N\setminus B_\varepsilon(x)} \dfrac{|u(x)-u(y)|^{p-2}(u(x)-u(y))}{|x-y|^{N+sp}}dy.$$
  One typical feature of problem \hyperref[P]{$(P_p)$} is the nonlocality, in the sense that the value of $(-\Delta)^s_pu(x)$ at any point
  $x\in \Omega$ depends not only on $\Omega$, but actually on the entire space $\R^N$. In recent years, the  problem \hyperref[P]{$(P_p)$} has been studied in many papers, we refer to \cite{amm,frr,lii}, in which the authors have used different methods to get the existence of solutions for \hyperref[P]{$(P_p)$}.
  
  This paper is organized as follows : in the second section, we recall some properties of Orlicz-Sobolev and fractional Sobolev spaces.  The third section is devoted to proving some properties of the
   fractional Orlicz-Sobolev spaces. 
  Finally,   using the direct method in calculus variations, we obtain the existence of a weak solution of problem  \hyperref[P_a]{$(P_a)$}.

\section{Some preliminaries results}$\label{100}$
First, we briefly recall the definitions and some elementary properties of  the Orlicz-Sobolev spaces. We refer the reader to \cite{1,28,34} for further reference and for some of the proofs of the results in this section.
\subsection{Orlicz-Sobolev Spaces} 
 We start by recalling the definition of the well-known N-functions. 
 
  Let $\Omega$ be an open subset of $\R^N$. Let $A$ :
  $\R^+ \rightarrow \R^+$ be an N-function, that is, $A$ is continuous, convex, with $A(t) > 0$ for $t >
  0$, $ \frac{A(t)}{t}\rightarrow 0$ as  $t \rightarrow 0$ and $ \frac{A(t)}{t} \rightarrow\infty$ as $t \rightarrow \infty$. Equivalently, $A$ admits the
  representation : $ A(t)=\int_{0}^{t} a(s)ds$ where $a : \R^+ \rightarrow \R^+$ is non-decreasing, right
  continuous, with $a(0) = 0$, $ a(t) > 0$ $\forall t> 0$ and $a(t) \rightarrow \infty $ as $t \rightarrow \infty$.  The conjugate N-function of $A$ is defined by $\overline{A}(t) = \int_{0}^{t} \overline{a}(s)ds$, where $ \overline{a} : \R^+\rightarrow
  \R^+$ is given by $\overline{a}(t) = \sup \left\lbrace s : a(s) \leqslant t\right\rbrace$.  Evidently we have
  \begin{equation}\label{1}
st\leqslant A(t)+\overline{A}(s),
  \end{equation}
  which is known Young's inequality. Equality holds in (\ref{1}) if and only if either $t=\overline{a}(s)$ or $s=a(t)$.\\
   We will extend these
  N-functions into even functions on all $\R$. The N-function $A$ is said to satisfy the
  global $\Delta_2$-condition if, for some $k > 0,$
 $$ A(2t) \leqslant kA(t)\text{ , }\forall t \geqslant 0.$$
  When this inequality holds only for $t \geqslant t_0 > 0$, $A$ is said to satisfy the $\Delta_2$-condition
  near infinity.
\\
We call the pair $(A,\Omega)$ is $\Delta$-regular if either :\\
(a) $A$  satisfies a
  global $\Delta_2$-condition, or \\
  (b) $A$ satisfies a $\Delta_2$-condition near infinity and $\Omega$ has finite volume.\\
  Throughout this paper, we assume that
    \begin{equation}\label{A}
      1<p_0:=\inf_{s> 0}\dfrac{sa(s)}{A(s)}<p^0:=\sup_{s> 0}\dfrac{sa(s)}{A(s)}<+\infty.
          \end{equation}
          which assures that $A$ satisfies the global $\Delta_2$-condition.
          
    \begin{lemma}$\label{2.1..}$ (see. \cite{9}).
      Let $A$ be an N-function which satisfies the global $\Delta_2$-condition. Then we have,
      \begin{equation}\label{2}
      \overline{A}(a(t))\leqslant c A(t) \text{ for all } t\geqslant 0
      \end{equation} 
      where $c>0$.
     \end{lemma}
   Let $\Omega$ be an open subset of $\R^N$. The Orlicz class $K_A (\Omega)$ (resp. the Orlicz space
   $L_A(\Omega)$) is defined as the set of (equivalence classes of) real-valued measurable
   functions $u$ on $\Omega$ such that  
\begin{equation}\label{3}
\int_{\Omega} A(|u(x)|)dx <\infty \hspace*{0.5cm} \text{ (resp. } \int_{\Omega} A(\lambda |u(x)|)dx< \infty \text{ for some } \lambda >0 ).
\end{equation}
$L_A(\Omega)$ is a Banach space under the Lexumburg norm
\begin{equation}\label{4}
 ||u||_A=\inf \Bigg\{\lambda>0  : \int_{\Omega}A\left(  \dfrac{|u(x)|}{\lambda}\right) dx \leqslant 1\Bigg\},
\end{equation}
and $K_A(\Omega)$ is a convex subset of $L_A(\Omega)$. The closure in $L_A(\Omega)$ of the set of
bounded measurable functions on $\Omega$ with compact support in $\overline{\Omega}$ is denoted by $E_A(\Omega)$.\\
The equality $E_A(\Omega)=L_A(\Omega)$ holds if and only if $(A,\Omega)$ is $\Delta$-regular.

 Using the Young's inequality, it is possible to prove a H\"older type inequality,
  that is,
  \begin{equation}
  \left| \int_{\Omega}uvdx\right| \leqslant 2||u||_A||v||_{\overline{A}}\hspace*{0.5cm} \text{ for all }  u \in L_A(\Omega)  \text{ and all } v\in L_{\overline{A}}(\Omega).
  \end{equation}
\subsection{Fractional Sobolev spaces} 
This subsection is devoted to the definition of the fractional Sobolev spaces, and we recall some result of continuous and compact embedding of fractional Sobolev spaces. We refer the reader to \cite{20,24} for further reference and for some of the proofs of these results.

We start by fixing the fractional exponent $s\in(0,1)$. For any $p\in [1,\infty)$, we define the fractional Sobolev space $W^{s,p}(\Omega)$ as follows :
$$ W^{s,p}(\Omega)=\Bigg\{u \in L^p(\Omega)  \text{ : }  \dfrac{|u(x)-u(y)|}{|x-y|^{\frac{N}{p}+s}} \in L^p(\Omega \times \Omega)\Bigg\},$$
that is, an intermediary Banach space, endowed
with its natural norm
$$||u||_{s,p}=\Bigg(\int_{\Omega}|u|^pdx+\int_{\Omega}\int_{\Omega} \dfrac{|u(x)-u(y)|^p}{|x-y|^{sp+N}}dxdy\Bigg)^{\frac{1}{p}}.$$
 \begin{theorem}$\label{2.3}$(see. \cite{24}).
   Let $s\in (0,1)$, $p\in [1,+\infty)$ and let $\Omega$  be an open
    subset of  $\R^N$  with $C^{0,1}$-regularity and bounded boundary. Then there exists a constant $C=C(N,s,p,\Omega)$ such that, for all  $f \in W^{s,p}(\Omega)$, we have  \\
     $$||f||_{L^q(\Omega)} \leqslant C||f||_{W^{s,p}(\Omega)} \text{  for all  }  q \in [p,p^*],$$
     that is,
     $$ W^{s,p}(\Omega)  \hookrightarrow L^q(\Omega)  \text{  for all   }  q \in[p,p^*],$$
       where  {\small$$
                p^*= \hspace*{0.1cm} \left\{ 
                  \begin{array}{clclc}
               \frac{Np}{N-sp}\hspace*{0.5cm} \text{ if } N>sp
                 \\\\
                   \infty\hspace*{0.5cm}  \text{ if }  N \leqslant sp.
                   \label{eq1}
                  \end{array}
                  \right. 
               $$ }\\
       If, in addition, $\Omega$ is bounded, then the space $W^{s,p}(\Omega)$ is continuously
       embedded in $L^q(\Omega)$ for any $q\in [1,p^*]$.
  \end{theorem}
  \begin{theorem}$\label{2.4}$(see. \cite{20}).
   Let $s\in (0,1)$, $p\in [1,+\infty)$ and let $\Omega$  be a bounded open
    subset of  $\R^N$   with $C^{0,1}$-regularity and bounded boundary. Then \\   
  (i) if $sp<N$, then the embedding $W^{s,p}(\Omega)  \hookrightarrow L^q(\Omega)$ is compact   for all $q \in[1,p^*)$;\\  
    (ii) if $sp=N$, then the embedding $W^{s,p}(\Omega)  \hookrightarrow L^q(\Omega)$ is compact   for all $q \in[1,\infty)$;\\    
        (iii) if $sp>N$, then the embedding $W^{s,p}(\Omega)  \hookrightarrow L^\infty(\Omega)$ is compact.
  \end{theorem}
  
 \begin{theorem}\label{AA2.2} (see. \cite{110})
     Suppose that $X$ is a reflexive Banach space with norm $||.||$ and let
     $V\subset X$ be a weakly closed subset of $X$. Suppose $E : V \longrightarrow \R \cup \left\lbrace +\infty\right\rbrace $ is coercive
     and (sequentially) weakly lower semi-continuous on $V$ with respect to $X$, that
     is, suppose the following conditions are fulfilled:
     
     (1) $E(u)\rightarrow \infty$ as $||u||\rightarrow \infty$, $u\in V$.
     
     (2)  For any $u\in V$, any sequence $\left\lbrace u_n\right\rbrace $ in $V$ such that $u_n\rightharpoonup u$ weakly in $X$
     there holds:
    $$E(u)\leqslant \liminf_{n\rightarrow \infty}E(u_n).$$
     Then $E$ is bounded from below on $V$ and attains its infinitum in $V$.     
      \end{theorem}
  \section{Variational framework}
       Now, we define the fractional Orlicz-Sobolev spaces, and we will present some important results of these spaces.
       \begin{definition}
       Let $A$ be an N-function. For a given domain $\Omega$ in $\R^N$ and $0<s<1$, we define  the fractional Orlicz-Sobolev space $W^sL_A(\Omega)$ as follows :
       \begingroup\makeatletter\def\f@size{10}\check@mathfonts\begin{equation}\label{5}
       W^s{L_A}(\Omega)=\Bigg\{u\in L_A(\Omega) :  \int_{\Omega} \int_{\Omega} A\left( \dfrac{\lambda| u(x)- u(y)|}{|x-y|^s}\right)\dfrac{ dxdy}{|x-y|^N}< \infty \text{ for some } \lambda >0 \Bigg\}.
       \end{equation}
       \endgroup
       This space is equipped with the norm,
       \begin{equation}\label{6}
       ||u||_{s,A}=||u||_{A}+[u]_{s,A},
       \end{equation}
       where $[.]_{s,A}$ is the Gagliardo seminorm, defined by 
       \begin{equation}\label{7}
       [u]_{s,A}=\inf \Bigg\{\lambda >0 :  \int_{\Omega} \int_{\Omega} A\left( \dfrac{| u(x)- u(y)|}{\lambda|x-y|^s}\right)\dfrac{ dxdy}{|x-y|^N}\leqslant 1 \Bigg\}.
       \end{equation}
       \end{definition}
       \begin{definition}
       Let $A$ be an N-function. For a given domain $\Omega$ in $\R^N$ and $0<s<1$, We define, the space $W^sE_A(\Omega)$ as follows :
       \begin{equation}\label{8}
       W^s{E_A}(\Omega)=\left\lbrace u\in E_A(\Omega) : \int_{\Omega} \int_{\Omega} A\left( \dfrac{| u(x)- u(y)|}{|x-y|^s}\right)\dfrac{ dxdy}{|x-y|^N}<\infty \right\rbrace.
       \end{equation}
       \end{definition}
       \begin{remark} \text{  }\\
       $\bullet$ $W^s{E_A}(\Omega)\subset W^s{L_A}(\Omega)$.\\
        $\bullet$ $W^s{E_A}(\Omega)$ coincides with  $W^s{L_A}(\Omega)$ if and only if  $(A,\Omega)$ is $\Delta$-regular.\\ 
        $\bullet$ If $1<p<\infty$ and $A_p(t)=t^p$, then $W^sL_{A_p}(\Omega)=W^sE_{A_p}(\Omega)=W^{s,p}(\Omega)$.
        \end{remark}
       
      To simplify the notation, we put           
      $$     h_{x,y}(u):=\dfrac{| u(x)- u(y)|}{|x-y|^s}.$$  
        \begin{theorem}
       Let $\Omega$ be an open subset of $\R^N$, and let $s\in (0,1)$. The space $W^sL_A(\Omega)$ is a Banach space with respect to the norm $(\ref{6})$, and a  separable (resp. reflexive) space if and only if $(A,\Omega)$ is $\Delta$-regular (resp. $(A,\Omega)$ and $(\overline{A},\Omega)$ are $\Delta$-regular).  Furthermore
       if   $(A,\Omega)$ is $\Delta$-regular and $A(\sqrt{t})$ is convex, then  the space $W^sL_A(\Omega)$ is uniformly convex.
       \end{theorem}
        \noindent \textit{Proof}. 
        Let $\left\lbrace u_n\right\rbrace $ be a Cauchy sequence for the norm $||.||_{s,A}$. In particular, $\left\lbrace u_n\right\rbrace $
        is a Cauchy sequence in $L_A(\Omega)$. It converges to a function $u\in L_A(\Omega)$. Moreover, the
        sequence $h_{x,y}(u_n)$
        is a Cauchy sequence in $L_A(\Omega\times \Omega,d\mu)$, where $\mu$   is a  measure on  $\Omega\times\Omega$ which is given by
        $$d\mu :=|x-y|^{-N}dxdy.$$
         It therefore also converges to an element of  $L_A(\Omega\times \Omega,d\mu)$.
        Let us extract a subsequence $\left\lbrace u_{\sigma(n)}\right\rbrace $ of $\left\lbrace u_n\right\rbrace $ that converges almost everywhere to $u$. We note that $h_{x,y}(u_{\sigma(n)})$ converges, for almost every pair $(x, y)$ to $h_{x,y}(u)$.
        Applying Fatou's lemma, we obtain, for some $\lambda$ (note that $\lambda$ exists since $\left\lbrace u_{\sigma(n)}\right\rbrace  \subset W^sL_A(\Omega)$),
                \begingroup\makeatletter\def\f@size{9}\check@mathfonts$$ \int_{\Omega} \int_{\Omega}  A\left( \dfrac{\lambda |u(x)- u(y)|}{|x-y|^s}\right) \dfrac{dxdy}{|x-y|^N}\leqslant \liminf_{n\rightarrow \infty} \int_{\Omega} \int_{\Omega} A\left( \dfrac{\lambda| u_{\varphi(n)}(x)- u_{\varphi(n)}(y)|}{|x-y|^s}\right)  \dfrac{dxdy}{|x-y|^N}<\infty.$$\endgroup        
       Hence $u\in W^s{L_A}(\Omega)$.\\
       On the other hand, since $h_{x,y}(u_n)$ converges in $L_A(\Omega\times \Omega,d\mu)$, then by dominated convergence theorem, there exist a subsequence $h_{x,y}(u_{\sigma(n)})$ and a function $ k$ in $L_A(\Omega\times \Omega,d\mu)$ such that 
       $$|h_{x,y}(u_{\sigma(n)})|\leqslant |k(x,y)| \text{ for almost every pair $(x,y)$,} $$
       and we have 
       $$h_{x,y}(u_{\sigma(n)})\longrightarrow h_{x,y}(u)  \text{ for almost every pair $(x,y)$,} $$
       this implies by dominated convergence theorem that, $$[u_n-u]_{s,A}\longrightarrow 0.$$
       Finally $u_n\rightarrow u$ in $W^sL_A(\Omega)$.
       
        To establish the reflexivity and separation of the fractional  Orlicz-Sobolev  spaces, we define the operator T : $W^sL_A(\Omega)\rightarrow L_{A}(\Omega) \times  L_A(\Omega \times \Omega,d\mu)$ by \\
        $$T(u)=\left( u(x), \dfrac{|u(x)-u(y)|}{|x-y|^s}\right) .$$
        Clearly, T is an isometry. Since $L_A(\Omega)$ is a reflexive, separable space and uniformly convex (see \cite{1,A10.}), then $W^sL_A(\Omega)$ is also a reflexive, separable space and uniformly convex.\hspace*{11cm $\Box$ } 
               
                           Let $\widetilde{W}^s_0L_A(\Omega)$
             denote the closure of  $C^{\infty}_0(\Omega)$ in the norm $||.||_{s,A}$ defined
            in $(\ref{6})$. Then we have the following result.
                       \begin{theorem}\label{33}(Generalized Poincar\'e inequality). 
                Let $\Omega$ be a bounded open subset of  $\R^N,$ and let $s\in (0,1)$. Let $A$ be an N-function. Then there exists a positive constant $\mu$ such that,        $$ ||u||_{A}\leqslant \mu [u]_{s,A} \text{ for all  }  u \in  \widetilde{W}^s_0L_A(\Omega).$$
                               \end{theorem}
                 Therefore, if $\Omega$ is bounded and $A$ be an N-function, then $ [.]_{s, A}$ is a norm of $\widetilde{W}^s_0L_A(\Omega) $ equivalent to $ ||. ||_{s,A}.$
                   \\              
                 \noindent \textit{Proof of Theorem \ref{33}}. 
                 Since $\widetilde{W}^s_0L_A(\Omega)$ is the closure of $C_0^{\infty}(\Omega)$ in $W^s{L_A}(\Omega)$, then it is enough to prove that there exists a positive constant $\mu$ such that, 
                 $$||u||_{A}\leqslant \mu [u]_{s,A} \text{ for all   }  u \in  C_0^{\infty}(\Omega).$$
                 Indeed, let $ u \in C_0^{\infty}(\Omega) $ and $ B_R \subset \R ^ N \setminus \Omega $, that is, the ball of radius $ R $ in the complement of $ \Omega $. Then for all $x\in \Omega$, $y\in B_R$ and all $\lambda>0$ we have,
                 $$A\left( \dfrac{|u(x)|}{\lambda}\right) =A\left( \dfrac{|u(x)-u(y)|}{\lambda |x-y|^s }|x-y|^s \right)\dfrac{|x-y|^N}{|x-y|^N},$$
                 this implies that,
                        \begingroup\makeatletter\def\f@size{10}\check@mathfonts  $$A\left( \dfrac{|u(x)|}{\lambda}\right) \leqslant A\left( \dfrac{|u(x)-u(y)|}{\lambda |x-y|^s} diam(\Omega\cup B_R)^s \right)\dfrac{diam(\Omega\cup B_R)^N}{|x-y|^N},$$\endgroup
       we suppose $\alpha =diam(\Omega\cup B_R)^s$, we get
       $$A\left( \dfrac{|u(x)|}{\alpha\lambda}\right) \leqslant A\left( \dfrac{|u(x)-u(y)|}{\lambda |x-y|^s} \right)\dfrac{diam(\Omega\cup B_R)^N}{|x-y|^N},$$
       therefore
               {\small $$|B_R|A\left( \dfrac{|u(x)|}{\alpha\lambda}\right) \leqslant diam(\Omega\cup B_R)^N\int_{B_R} A\left( \dfrac{|u(x)-u(y)|}{\lambda |x-y|^s} \right) \dfrac{dy}{|x-y|^N},$$}
       then
                  $$ \int_{\Omega} A\left( \dfrac{|u(x)|}{\alpha\lambda}\right) dx \leqslant \dfrac{diam(\Omega\cup B_R)^N}{|B_R|}\int_{\Omega}\int_{B_R} A\left( \dfrac{|u(x)-u(y)|}{\lambda |x-y|^s} \right)\dfrac{dxdy}{|x-y|^N},$$ 
                  so, 
                 $$ ||u||_{A}\leqslant \mu [u]_{s,A} \text{ for all  }  u \in  C^{\infty}_0(\Omega),$$
                 where $\mu =\dfrac{ diam(\Omega\cup B_R)^N \alpha}{|B_R|}$. 
                 By passing to the limit, the desired result is obtained.
                      \hspace*{12.7cm$\Box$ }  
                      \begin{corollary}
                       Let $\Omega$ be a bounded open subset of  $\R^N,$ and let $s\in (0,1)$. Let $A$ be an N-function. We define the space $W^s_0L_A(\Omega)$ as follows :
                       $$W^s_0L_A(\Omega)=\left\lbrace u\in W^sL_A(\R^N) \text{ : } u=0 \text{ a.e in } \R^N\setminus \Omega \right\rbrace, $$             
                        Then there exists a positive constant $\mu$ such that,          $$ ||u||_{A}\leqslant \mu [u]_{s,A} \text{ for all   }  u \in  W^s_0L_A(\Omega).$$
                      \end{corollary}  
     Proof of this corollary is similar to proof of Theorem $\ref{33}$.
                 \begin{theorem}      \label{3.1..}
            Let $\Omega$ be a bounded open subset of $\R^N$, let $0<s<1$ and let $A$ be an $N$-function. 
            Then,
            $$C^2_0(\Omega)\subset W^s_0L_A(\Omega).$$
       \end{theorem}
       \begin{lemma}\label{lem1}
       Let $A$ be an $N$-function. Then 
       $$\dfrac{\delta(x)}{|x|^s}\in L_A(\R^N,|x|^{-N}dx) \text{ with }    \delta(x)=\min\left\lbrace 1,|x|\right\rbrace.$$
       \end{lemma}
                  \noindent \textit{Proof of Lemma \ref{lem1}}. 
         We put $\Omega_1=\left\lbrace x\in \R^N \text{ : } |x|>1\right\rbrace $  and    $\Omega_2=\left\lbrace x\in \R^N \text{ : } |x|\leqslant 1\right\rbrace $. Then, we have
          $$
                \begin{aligned}
               \int_{\R^N}A\left( \dfrac{\delta(x)}{|x|^s}\right)  \dfrac{dx}{|x|^N} & = \int_{\Omega_1}A\left( \dfrac{\delta(x)}{|x|^s}\right)  \dfrac{dx}{|x|^N}+\int_{\Omega_2}A\left( \dfrac{\delta(x)}{|x|^s} \right) \dfrac{dx}{|x|^N}\\
               & =\int_{\Omega_1}A\left( \dfrac{1}{|x|^s}\right)  \dfrac{dx}{|x|^N}+\int_{\Omega_2}A\left( \dfrac{|x|}{|x|^s} \right) \dfrac{dx}{|x|^N}\\
               & \leqslant A\left( 1 \right)\int_{\Omega_1}  \dfrac{dx}{|x|^{N+s}}+A\left( 1 \right)\int_{\Omega_2} \dfrac{dx}{|x|^{N+s-1}},\\
                  \end{aligned} 
                       $$
              note that the last integrals are finite since $N+s>N$ and $N+s-1<N$ respectively. Therefore 
              $$\int_{\R^N}A\left( \dfrac{\delta(x)}{|x|^s}\right)  \dfrac{dx}{|x|^N}<\infty.$$
                       
                  \noindent \textit{Proof of Theorem \ref{3.1..}}. 
            Let $u\in C^2_0(\Omega)$, we only need that to check that 
            $$ \int_{\R^N} \int_{\R^N} A\left( \dfrac{\lambda| u(x)- u(y)|}{|x-y|^s}\right) \dfrac{dxdy}{|x-y|^N} <\infty \text{ for some } \lambda>0.$$
            Indeed, since $u$ vanishes outside $\Omega$, we have
             \begingroup\makeatletter\def\f@size{10.5}\check@mathfonts $$
             \begin{aligned}
            \int_{\R^N}\int_{\R^N}   A\left( \dfrac{| u(x)- u(y)|}{|x-y|^s}\right) \dfrac{dxdy}{|x-y|^N} &= \int_{\Omega} \int_{\Omega} A\left( \dfrac{| u(x)- u(y)|}{|x-y|^s}\right) \dfrac{dxdy}{|x-y|^N}\\ &\hspace*{0.3cm}+2\int_{\Omega}\int_{\R^N \setminus \Omega}   A\left( \dfrac{| u(x)- u(y)|}{|x-y|^s}\right)\dfrac{dxdy}{|x-y|^N}\\
            &\leqslant 2\int_{\Omega}\int_{\R^N}   A\left( \dfrac{| u(x)- u(y)|}{|x-y|^s}\right) \dfrac{dxdy}{|x-y|^N}.
            \end{aligned} 
            $$\endgroup
            Now we notice that
            $$|u(x)-u(y)|\leqslant ||\nabla u||_{L^{\infty}(\R^N)}|x-y| \text{ and } |u(x)-u(y)|\leqslant 2||u||_{L^{\infty}(\R^N)}.$$
            Accordingly, we get 
            $$|u(x)-u(y)|\leqslant 2||u||_{C^{1}(\R^N)}\min\left\lbrace 1,|x-y|\right\rbrace:=\alpha \delta(x-y), $$
            with $\alpha =2||u||_{C^{1}(\R^N)}$ and since $\dfrac{\delta(x)}{|x|^s}\in L_A(\R^N,|x|^{-N}dx)$. There exists $\lambda>0$, such that,
             \begingroup\makeatletter\def\f@size{10}\check@mathfonts $$
             \begin{aligned}
             \int_{\R^N}\int_{\R^N}   A\left( \dfrac{\lambda| u(x)- u(y)|}{\alpha|x-y|^s}\right) \dfrac{dxdy}{|x-y|^N} & \leqslant 2\int_{\Omega}\int_{\R^N}  A\left( \lambda \dfrac{\delta(x-y)}{|x-y|^s}\right) \dfrac{dxdy}{|x-y|^N}\\
             & = 2|\Omega|\int_{\R^N}  A\left( \lambda \dfrac{\delta(\xi)}{|\xi|^s}\right) \dfrac{d\xi}{|\xi|^N}<\infty,
             \end{aligned}
             $$\endgroup
            this implies that  $u\in W^s_0L_A(\Omega)$.
                                \hspace*{8.4cm$\Box$ }  
\begin{remark}
A trivial consequence of Theorem $\ref{3.1..}$, $W^sL_A(\Omega)$  is non-empty.
\end{remark}        
             \begin{proposition}\label{pro3}
                    Let $\Omega$ be an open subset of $\R^N$ and let $A$ be an $N$-function. Assume condition $(\ref{A}) $ is satisfied, then the following relations hold true 
                    
                    \begin{equation}\label{32}
                    [u]^{p_0}_{s,A}\leqslant \phi(u) \leqslant [u]^{p^0}_{s,A} \text{    for all  }  u\in W^sL_A(\Omega) \text{ with }[u]_{s,A}>1,
                    \end{equation} 
                     \begin{equation}\label{A33}
                       [u]^{p^0}_{s,A}\leqslant \phi(u) \leqslant [u]^{p_0}_{s,A} \text{   for all  }  u\in W^sL_A(\Omega) \text{ with }[u]_{s,A}<1,
                       \end{equation} 
         where $\phi (u)=  $$\displaystyle\int_{\Omega} \int_{\Omega} A\left( \dfrac{| u(x)- u(y)|}{|x-y|^s}\right)\dfrac{dxdy}{|x-y|^N}$.       
                    \end{proposition}
                    \noindent \textit{Proof}.
                     First we show that  $\phi(u) \leqslant [u]^{p^0}_{s,A} \text{   for all  }  u\in W^sL_A(\Omega) \text{ with }[u]_{s,A}>1$. Indeed, since $p^0\geqslant \dfrac{ta(t)}{A(t)} $ for all $t>0$ it follows that for letting $\sigma >1$ we have 
                    $$ \log(A(\sigma t))-\log(A(t))=\int_{t}^{\sigma t} \dfrac{a(\tau)}{A(\tau)}d\tau   \leqslant \int_{t}^{\sigma t} \dfrac{p^0}{\tau}d\tau=\log(\sigma^{p^0}).$$ 
                    Thus, we deduce 
                    \begin{equation}\label{34}
                    A(\sigma t)\leqslant \sigma^{p^0}A(t) \text{ for all } t>0 \text{ and } \sigma>1.
                    \end{equation}
                    Let now $u\in W^sL_A(\Omega)$ with $[u]_{s,A}>1$. Using the definition of the Luxemburg norm and the relation $(\ref{34})$, we deduce 
                    $$
                    \begin{aligned}
                    \int_{\Omega}\int_{\Omega}A\left( h_{x,y}(u)\right) \dfrac{dxdy}{|x-y|^N}&=\int_{\Omega}\int_{\Omega}A\left( [u]_{s,A} \dfrac{h_{x,y}(u)}{[u]_{s,A}}\right) \dfrac{dxdy}{|x-y|^N}\\
                    &\leqslant [u]^{p^0}_{s,A}\displaystyle\int_{\Omega}\int_{\Omega}A\left( \dfrac{h_{x,y}(u)}{[u]_{s,A}} \right) \dfrac{dxdy}{|x-y|^N}\\
                    &\leqslant [u]^{p^0}_{s,A}.
                    \end{aligned}
                    $$
                    Now, we show that $\phi(u) \geqslant [u]^{p_0}_{s,A} \text{   for all  } u\in W^sL_A(\Omega) \text{ with }[u]_{s,A}>1$. Indeed since $$p_0\leqslant \dfrac{ta(t)}{A(t)}$$ for all $t> 0$, it follows that letting $\sigma>1$ we have
                     $$ \log(A(\sigma t))-\log(A(t))=\int_{t}^{\sigma t} \dfrac{a(\tau)}{A(\tau)}d\tau   \geqslant \int_{t}^{\sigma t} \dfrac{p_0}{\tau}d\tau=\log(\sigma^{p_0}).$$ 
                      Thus, we deduce 
                        \begin{equation}\label{35}
                        A(\sigma t)\geqslant \sigma^{p_0}A(t) \text{ for all } t>0 \text{ and } \sigma>1.
                        \end{equation}
                   Let $u\in W^sL_A(\Omega)$ with $[u]_{s,A}>1$, we consider $\beta\in (1,[u]_{s,A})$, since $\beta<[u]_{s,A}$, so by definition of Luxemburg norm, it follows that 
                 $$\int_{\Omega}\int_{\Omega}A\left( \dfrac{h_{x,y}(u)}{\beta} \right) \dfrac{dxdy}{|x-y|^N}>1,$$
                 the above consideration implies that 
                 $$
                 \begin{aligned}
                 \displaystyle\int_{\Omega}\int_{\Omega}A\left( h_{x,y}(u)\right) \dfrac{dxdy}{|x-y|^N} &=\int_{\Omega}\int_{\Omega}A\left( \beta\dfrac{h_{x,y}(u)}{\beta} \right) \dfrac{dxdy}{|x-y|^N}\\
                 &\geqslant \beta^{p_0}\int_{\Omega}\int_{\Omega}A\left( \dfrac{h_{x,y}(u)}{\beta} \right) \dfrac{dxdy}{|x-y|^N}\\
                 &\geqslant \beta^{p_0},
                 \end{aligned}
                 $$
                 letting $\beta\nearrow [u]_{s,A}$, we deduce that relation $(\ref{32})$ hold true.
                 
                 Next, we show that $\phi(u) \leqslant [u]^{p_0}_{s,A} \text{   for all  }  u\in W^sL_A(\Omega) \text{ with }[u]_{s,A}<1$. Similar techniques as those used in the proof of relation $(\ref{34})$ and $(\ref{35})$, we have 
                 \begin{equation}\label{36}
                 A(t)\leqslant \tau^{p_0}A\left( \dfrac{t}{\tau} \right) \text{ for all } t>0 \text{ and } \tau \in(0,1).
                 \end{equation}
                  Let $u\in W^sL_A(\Omega)$ with $[u]_{s,A}<1$.
                  Using the definition of the Luxemburg norm and the relation $(\ref{36})$, we deduce 
                     $$
                     \begin{aligned}
                     \displaystyle\int_{\Omega}\int_{\Omega} A\left( h_{x,y}(u)\right) \dfrac{dxdy}{|x-y|^N}
                     &\leqslant [u]^{p_0}_{s,A}\displaystyle\int_{\Omega}\int_{\Omega}A\left( \dfrac{h_{x,y}(u)}{[u]_{s,A}} \right) \dfrac{dxdy}{|x-y|^N}\\
                     &\leqslant [u]^{p_0}_{s,A}.
                     \end{aligned}
                     $$
                     
                     Finally, we show that $\phi(u) \geqslant [u]^{p^0}_{s,A} \text{   for all  }  u\in W^sL_A(\Omega) \text{ with }[u]_{s,A}<1$. Similar techniques as those used in the proof of relation $(\ref{34})$ and $(\ref{35})$, we have 
                     \begin{equation}\label{37}
                     A(t)\geqslant \tau^{p^0}A\left( \dfrac{t}{\tau} \right) \text{ for all } t>0 \text{ and } \tau \in(0,1).
                     \end{equation}
                      Let $u\in W^sL_A(\Omega)$ with $[u]_{s,A}<1$ and $\beta\in (0,[u]_{s,A})$, so by $ (\ref{37})$ we find 
                 \begin{equation}\label{38}
                 \displaystyle\int_{\Omega}\int_{\Omega}A\left(h_{x,y}(u)\right) \dfrac{dxdy}{|x-y|^N}\geqslant \beta^{p^0}\int_{\Omega}\int_{\Omega}A\left( \dfrac{h_{x,y}(u)}{\beta} \right) \dfrac{dxdy}{|x-y|^N}.
                 \end{equation}
                      We define $v(x)=\dfrac{u(x)}{\beta}$ for all $x\in \Omega$, we have $[v]_{s,A}=\dfrac{[u]_{s,A}}{\beta}>1$. Using the relation $(\ref{32})$ we find
                      \begin{equation}\label{39}
                    \int_{\Omega}\int_{\Omega}A\left( \dfrac{h_{x,y}(u)}{\beta} \right) \dfrac{dxdy}{|x-y|^N}=  \displaystyle\int_{\Omega}\int_{\Omega}A\left( h_{x,y}(v)\right) \dfrac{dxdy}{|x-y|^N}>[v]^{p_0}_{s,A}>1,
                      \end{equation}
                      by $(\ref{38})$ and $(\ref{39})$ we obtain 
                      $$\displaystyle\int_{\Omega}\int_{\Omega}A\left( h_{x,y}(u)\right) \dfrac{dxdy}{|x-y|^N}\geqslant \beta^{p^0}.$$
                      Letting $\beta\nearrow [u]_{s,A}$, we deduce that relation $(\ref{A33})$ hold true.
                     \hspace*{3.6cm$\Box$ }                   
                     
       The embeddings results obtained in the fractional Sobolev space $W^{s,p}(\Omega)$ can also be formulated for the fractional Orlicz-Sobolev spaces.

         \begin{theorem}\label{3.4.}
            Let $\Omega$ be a bounded open subset of $\R^N$.   Let $\Omega$ be an open subset of $\R^N$ and let $A$ be an $N$-function and let $s\in (0,1)$. Assume condition $(\ref{A}) $ is satisfied, then the space $W^sL_{A}(\Omega)$ is continuously embedded in $W^{s,p_0}(\Omega)$.
            \end{theorem}
             \noindent \textit{Proof}.    By (\ref{A}), we have
           $$p_0:=\inf_{t> 0}\dfrac{ta(t)}{A(t)}>1,$$
           this implies that 
           $$\dfrac{p_0}{t}\leqslant \dfrac{a(t)}{A(t)},$$
           for all $t> 0$, it follows that letting $0<t_0<t$ we have
           $$\log(t^{p_0})-\log(t_0^{p_0})=\int_{t_0}^{ t} \dfrac{p_0}{\tau}d\tau \leqslant\int_{t_0}^{ t} \dfrac{a(\tau)}{A(\tau)}d\tau= \log(A( t))-\log(A(t_0)).$$             
           Therefore
           \begin{equation}\label{A24}
           |t|^{p_0}\leqslant cA(t) \hspace*{1cm} \text{ for all } t> 0.
           \end{equation}
           Let $u\in W^sL_A(\Omega)$ we have 
           $$\int_{\Omega} |u|^{p_0}dx\leqslant c\int_{\Omega}A(u)dx, $$
           then 
           \begin{equation}\label{A26}
           ||u||_{L^{p_0}}\leqslant c ||u||_{A},
           \end{equation}
           where $c>0$ and it is possibly different step by step.
           On the other hand, by estimation $(\ref{A24})$, we have 
            {\small$$
             \int_{\Omega}\int_{\Omega} \dfrac{|u(x)-u(y)|^{p_0}}{|x-y|^{p_0s+N}}dxdy\leqslant   c\int_{\Omega} \int_{\Omega} A\left( \dfrac{| u(x)- u(y)|}{|x-y|^s}\right)\dfrac{dxdy}{|x-y|^N}.
           $$}
            Then 
           \begin{equation}\label{A29}
           [u]_{s,p_0}\leqslant c [u]_{s,A}.
           \end{equation}
           Combining $(\ref{A26})$ and $(\ref{A29})$, we have
           $$ ||u||_{s,p_0}\leqslant c ||u||_{s,A}.$$
             The proof
             of Theorem \ref{3.4.} is complete.
                 \hspace*{7.1cm $\Box$ }\\
             By combining  Theorem \ref{3.4.} and Theorem \ref{2.4}, we obtain the following results.
           \begin{corollary}\label{3.1.}
           Let $\Omega$ be a bounded open subset of $\R^N$ with $C^{0,1}$-regularity and bounded boundary. Let $s\in (0,1)$ and let $A$ be an $N$-function satisfies the global $\Delta_2$-condition,
           we define  
           $$
        p_0^*= \hspace*{0.1cm} \left\{ 
                         \begin{array}{clclc}
                      \frac{Np_0}{N-sp_0}\hspace*{0.5cm} \text{ if } N>sp_0
                        \\\\
                          \infty\hspace*{0.5cm}  \text{ if }  N \leqslant sp_0.
                          \label{eq1}
                         \end{array}
                         \right. 
                   $$      
                  $\bullet$ If $sp_0<N$, then 
                  $W^sL_A(\Omega) \hookrightarrow L^{q}(\Omega),$ for all $q\in [1,p^*_0]$
                  and the embedding 
               $W^sL_A(\Omega) \hookrightarrow L^{q}(\Omega)$
               is compact for all $q\in [1 , p_0^*)$.\\      
              $\bullet$ If $sp_0=N$,  then 
                          $W^sL_A(\Omega) \hookrightarrow L^{q}(\Omega),$ for all $q\in [1,\infty]$
                          and the embedding 
                       $W^sL_A(\Omega) \hookrightarrow L^{q}(\Omega)$
                       is compact for all $q\in [1 , \infty)$.\\                
                       $\bullet$  If $sp_0>N$,  then the embedding 
    $W^sL_A(\Omega) \hookrightarrow L^{\infty}(\Omega)$
                                         is compact. 
                           \end{corollary} \section{Mains results}
   In this  section,  we prove the existence of a weak  solution for the problem \hyperref[P_a]{$(P_a)$} in fractional Orlicz Sobolev spaces, by means of  the direct method in calculus of variations. For this, we suppose
   that $f : \Omega \times \R \rightarrow \R$ is a Carath\'eodory function satisfying the following conditions :\\
                 there exist $\theta_1,\theta_2>0$, $1 < q < p_0^*$ and an open bounded set $\Omega_0\subset \Omega$  such that
                \label{f1} $$ (f_1) \hspace{4cm}|f(x,t)|\leqslant \theta_1 (1+|t|^{q-1}) \text{  a.e.  } (x,t)\in \Omega\times \R^N, \hspace{5cm} $$                     
                \label{f2} $$(f_2) \hspace{4cm}|f(x,t)|\geqslant \theta_2 |t|^{q-1} \text{  a.e.  } (x,t)\in \Omega_0\times \R^N. \hspace{5cm} $$
   
   In this section, we work in closed space $W^s_0L_A(\Omega)$ which can be equivalently renormed by setting $||.||:=[.]_{s,A}$.        
           \begin{definition}
          We say that $u\in W^s_0L_A(\Omega)$ is a weak solution of problem \hyperref[P_a]{$(P_a)$} if 
          $$<(-\Delta)^s_au,v> = \int_{\Omega} \int_{\Omega} a\left( |h_{x,y}(u)|\right)\dfrac{u(x)-u(y)}{|u(x)-u(y)|}h_{x,y}(v) \dfrac{dxdy}{|x-y|^{N} } =\int_{\Omega}f(x,u)vdx,$$
          for all $v\in W^s_0L_A(\Omega)$.
          \end{definition}
          \begin{theorem}\label{4.1.}
          Let  $A$  be an $N$- function satisfies $(\ref{A})$. Suppose that $f$ satisfies $\hyperref[f1]{(f_1)}$ and $\hyperref[f2]{(f_2)}$. If $1<q<p_0$, then the problem \hyperref[P_a]{$(P_a)$} has a nontrivial weak solution in $W^s_0L_A(\Omega)$. 
          \end{theorem}
          \begin{corollary}\label{cr}
          Let  $A$  be an $N$-function satisfies $(\ref{A})$. Suppose that $f$ satisfies $\hyperref[f1]{(f_1)}$ and $\hyperref[f2]{(f_2)}$. If $q=p_0$ and $\theta_1<\lambda_1/2$, where $\lambda_1$ is  define by
         \begin{equation}\label{lam}
          \lambda_1=\inf\limits_{u\in W^s_0L_A(\Omega)\setminus \left\lbrace 0\right\rbrace }\dfrac{||u||^{p_0}}{||u||_{p_0}^{p_0}}, \end{equation}
           then the problem \hyperref[P_a]{$(P_a)$} admits a weak solution in $W^s_0L_A(\Omega)$. 
          \end{corollary}
          For $u\in W^s_0L_A(\Omega)$, we define 
          $$J(u)= \int_{\Omega}\int_{\Omega}A\left( \dfrac{|u(x)-u(y)|}{|x-y|^s }\right)\dfrac{dxdy}{|x-y|^N},$$
          $$ H(u)=\int_{\Omega}F(x,u)dx \text{ and } I(u)=J(u)-H(u),$$
          where  $F(x,t)=\displaystyle\int_{0}^{t}f(x,t)d\tau$.
          Obviously the energy functional $I : W^s_0L_A(\Omega)\longrightarrow \R$ associated with problem \hyperref[P_a]{$(P_a)$} is well defined.
          \begin{lemma}\label{4.1}
          If $f$ satisfies assumption $\hyperref[f1]{(f_1)}$. Then the functional $H\in C^1(W^s_0L_A(\Omega),\R)$ and 
          \begin{equation}\label{fff}        
          <H'(u),v>=\int_{\Omega}f(x,u)vdx \text{ for all } u,v \in W^s_0L_A(\Omega).\end{equation}
          \end{lemma}
          \noindent \textit{Proof}. First, observe that by $(f_1)$ and the embedding theorem, $H$ is well-defined on $W_0^sL_A(\Omega)$. Usual arguments show  that $H$
           is G\^ateaux-differentiable on $W_0^sL_A(\Omega)$ with the derivative is given by $(\ref{fff})$.  Actually, let $\left\lbrace u_n\right\rbrace \subset W_0^sL_A(\Omega)$ be a sequence converging strongly to $u\in W_0^sL_A(\Omega)$. Since $W_0^sL_A(\Omega)$ is compactly embedded in $L^q(\Omega)$, then $\left\lbrace u_n\right\rbrace$ converges strongly to $u$ in $L^q(\Omega)$. So there exist a subsequence of $\left\lbrace u_n\right\rbrace$, still denoted by $\left\lbrace u_n\right\rbrace$, and a function $\overline{u} \in L^q(\Omega)$ such that $\left\lbrace u_n\right\rbrace$  converges to $u$ almost
           everywhere in $\Omega$ and $|u_n|\leqslant |\overline{u}|$ for all $n\in \mathbb{N}$ and almost everywhere in $\Omega$. By $(f_1)$, we have   for all measurable functions $u : \Omega\longrightarrow \mathbb{R}$, the operator defined by $u\longmapsto f(.,u(.))$  maps $L^q(\Omega)$ into $L^{q'}(\Omega)$. Fix $v\in W_0^sL_A(\Omega)$ with $||v||\leqslant 1$. By using the H\"older inequality and the embedding of $W_0^sL_A(\Omega)$ into $L^q(\Omega)$, we have 
           
           $$
           \begin{aligned}  
           |\left\langle  J'(u_n)-J'(u),v\right\rangle |& =\Bigg |\int_{\Omega}\left( f(x,u_n(x))-f(x,u(x))\right)v(x)dx \Bigg|,\\
           & \leqslant ||f(x,u_n(x))-f(x,u(x))||_{q'}||v||_{q},\\
           & \leqslant c_1||f(x,u_n(x))-f(x,u(x))||_{q'}||v||,
             \end{aligned}
             $$
             for some $c_1>0$. Thus, passing to the supremum for $||v||\leqslant 1$, we get 
             $$||J'(u_n)-J'(u)||_{(W^{s,p}_0(\Omega))^*}\leqslant c_1||f(.,u_n(.))-f(.,u(.))||_{q'}.$$
             Since $f$ is a continuous function, then  
             $$f(.,u_n(.))-f(.,u(.)) \longrightarrow 0 \text{  as  } n\rightarrow \infty$$
              and 
             $$|f(x,u_n(x))-f(x,u(x))|\leqslant \theta_1(2+|\overline{u}(x)|^{q-1}+|u(x)|^{q-1}).$$
            for almost everywhere $x\in \Omega$. Note that the majorant function in the previous relation is in $L^{q'}(\Omega)$. Hence, by
             applying the dominate convergence theorem we get that $||f(x,u_n(x))-f(x,u(x))||_{q'}\rightarrow 0$ as $n\rightarrow \infty$. This proves that $J'$  is continuous. 
                   \hspace*{13.3cm $\Box$ }
                           \begin{lemma}\label{4.2}
           The function $J\in C^1(W^s_0L_A(\Omega),\R)$ and
          $$<J'(u),v>= \int_{\Omega} \int_{\Omega} a\left( |h_{x,y}(u)|\right)\dfrac{u(x)-u(y)}{|u(x)-u(y)|}h_{x,y}(v) \dfrac{dxdy}{|x-y|^{N} }$$
          for all $u,v \in W^s_0L_A(\Omega)$. Moreover, for   each $u \in W^s_0L_A(\Omega)$, $J'(u) \in (W^s_0L_A(\Omega))^*$.

                  \end{lemma}
                  \noindent \textit{Proof}.
             First, it is easy to see that 
          \begin{equation}\label{30}
     <J'(u),v>= \int_{\Omega} \int_{\Omega} a\left( |h_{x,y}(u)|\right)\dfrac{u(x)-u(y)}{|u(x)-u(y)|}h_{x,y}(v) \dfrac{dxdy}{|x-y|^{N} }
          \end{equation}
           for all $u,v \in W^s_0L_A(\Omega)$. It follows from $(\ref{30})$  that $J'(u) \in (W^s_0L_A(\Omega))^*$ for each $u \in W^s_0L_A(\Omega)$.
           
           Next, we prove that $J\in C^1(W^s_0L_A(\Omega),\R)$. Let $\left\lbrace u_n\right\rbrace \subset W^s_0L_A(\Omega)$ with $u_n\longrightarrow u$ strongly in $W^s_0L_A(\Omega)$, then $h_{x,y}(u_n)\longrightarrow h_{x,y}(u)$ in $ L_{M}(\Omega\times\Omega,d\mu)$, so by dominated convergence theorem, there exist a subsequence $h_{x,y}(u_{n_k}) $ and a function $h$ in $ L_{A}(\Omega\times\Omega,d\mu)$ such that 
                                 $$
                    \begin{aligned}
                    U_{x,y}(u_{n_k}):=& a(|h_{x,y}(u_{n_k})|)\dfrac{u_{n_k}(x)-u_{n_k}(y)}{|u_{n_k}(x)-u_{n_k}(y)|}\\
                    & \longrightarrow  U_{x,y}(u):=a(|h_{x,y}(u)|) \dfrac{u(x)-u(y)}{|u(x)-u(y)|}  
                                   \end{aligned}
                    $$   
         and
        $$| U_{x,y}(u_{n_k})|\leqslant |a(|h|)|,$$
      for almost every $(x,y)$ in $\Omega\times\Omega$. By $(\ref{2})$,  $|a(h)| \in  L_{\overline{A}}(\Omega\times\Omega,d\mu)$, so for $v\in W^s_0L_A(\Omega)$  we have $h_{x,y}(v)\in L_{A}(\Omega\times\Omega,d\mu)$ 
                           and by H\"older's inequality,
                           \begingroup\makeatletter\def\f@size{10}\check@mathfonts       
                    $$     
                \begin{aligned}
               \left| \int_{\Omega}\int_{\Omega} (  U_{x,y}(u_{n_k})-
                U_{x,y}(u)) h_{x,y}(v)\dfrac{dxdy}{|x-y|^N}\right|&
               \leqslant 2\left|\left| U_{x,y}(u_{n_k})-
                          U_{x,y}(u)\right| \right|_{L_{\overline{A}}} \left| \left| h_{x,y}(v)\right| \right|_{L_{A}}\\
                 &= 2 \left|\left| U_{x,y}(u_{n_k})-
                            U_{x,y}(u)\right| \right|_{L_{\overline{A}}} \left| \left|v\right| \right|.
             \end{aligned}  
             $$ \endgroup
              Then by dominated convergence theorem we obtain that
              $$\sup_{||v||\leqslant 1} \Bigg| \int_{\Omega} \int_{\Omega}( U_{x,y}(u_{n_k})-
                         U_{x,y}(u)) h_{x,y}(v)\dfrac{dxdy}{|x-y|^N}\Bigg|\longrightarrow 0.$$  The proof of Lemma $\ref{4.2}$, is complete\hspace*{7.8cm $\Box$}\\ 
                               
           Combining Lemma $\ref{4.1}$ and Lemma $\ref{4.2}$, we get $I \in C^1(W^s_0L_A(\Omega),\R)$ and 
           {\small$$<I'(u),v>= \int_{\Omega} \int_{\Omega} a\left( |h_{x,y}(u)|\right)\dfrac{u(x)-u(y)}{|u(x)-u(y)|}h_{x,y}(v) \dfrac{dxdy}{|x-y|^{N} }- \int_{\Omega}f(x,u)vdx$$}
             for all  $u,v \in W^s_0L_A(\Omega).$     
             \begin{lemma}\label{4.3}
             Let  $\hyperref[f1]{(f_1)}$ be satisfied.  Then the functional $I\in C^1(W^s_0L_A(\Omega), \R)$ is weakly lower semicontinuous.
             \end{lemma}
            \noindent\textit{Proof}. First, note that the map :
             $$ u\longmapsto \int_{\Omega}\int_{\Omega}A\left( \dfrac{|u(x)-u(y)|}{|x-y|^s }\right)  \dfrac{dxdy}{|x-y|^{N} } $$ 
             is lower semicontinuous for the weak topology of $W^s_0L_A(\Omega)$. Indeed, 
             by Lemma $\ref{4.2}$, we have $J \in C^1(W^s_0L_A(\Omega), \R^N)$ and
             $$< J'(u),v>=\int_{\Omega} \int_{\Omega} a\left( |h_{x,y}(u)|\right)\dfrac{u(x)-u(y)}{|u(x)-u(y)|}h_{x,y}(v) \dfrac{dxdy}{|x-y|^{N} }$$
                for all  $u,v \in W^s_0L_A(\Omega)$. 
                 On the other hand, since $A$ is convex so $J$ is also convex.
                	Now, let $\left\lbrace u_n\right\rbrace \subset W^s_0L_A(\Omega)$ with $u_n\rightharpoonup u$ weakly in $ W^s_0L_A(\Omega)$. Then by convexity of $J$, we have 
                	$$J(u_n)-J(u)\geqslant <J'(u),u_n-u>,$$ 
                and	hence, we obtain $J(u)\leqslant \liminf J(u_n)$, that is, the map
           $$  u\longmapsto \int_{\Omega}\int_{\Omega}A\left( \dfrac{|u(x)-u(y)|}{|x-y|^s }\right)  \dfrac{dxdy}{|x-y|^{N} } $$      	
             is lower semicontinuous.
             
            Let $u_n\rightharpoonup u$ weakly in $W^s_0L_A(\Omega)$, so by Corollary $\ref{3.1.}$, $u_n\longrightarrow u$ in $L^q(\Omega)$ for all $q\in (1,p_0^*)$. Without loss of generality, we assume that $u_n\longrightarrow u$ a.e. in $\Omega$. Assumption $\hyperref[f1]{(f_1)}$ implies that 
            $$F(x,t)\leqslant \theta_1 (|t|^q+1).$$
            Thus, for any measurable subset $U\subset \Omega$,
            $$\int_{U}|F(x,u_n)|dx\leqslant \theta_1 \int_{U}|u_n|^qdx+\theta_1|U|.$$
            By H\"older inequality and corollary $\ref{3.1.}$ , we have,  
            \begin{equation}\label{37.}
                   \begin{aligned}
                  \int_{U}|F(x,u_n)|dx&\leqslant \theta_1 ||u_n^q||_{L^{\frac{p^*_0}{q}}}||1||_{L^{\frac{p^*_0}{p_0^*-q}}}+\theta_1|U|\\
                  &\leqslant \theta_1C ||u_n||^q|U|^{\frac{p_0^*-q}{p_0^*}}+\theta_1|U|.
                \end{aligned}  
            \end{equation}          
            It follows from $(\ref{37.})$ that the sequence $\left\lbrace |F(x,u_n)-F(x,u)|\right\rbrace $ is uniformly bounded and equi-integrable
            in $L^1(\Omega)$. The Vitali Convergence Theorem (see \cite{Ru}) implies
            $$\lim\limits_{n\rightarrow n}\int_{\Omega}|F(x,u_n)-F(x,u)|dx=0,$$
             so
             $$\lim\limits_{n\rightarrow \infty}\int_{\Omega}F(x,u_n)dx=\int_{\Omega}F(x,u)dx.$$
             Thus, the functional $H$ is weakly continuous. Further, we get that $I$ is weakly lower semicontinuous.
               \hspace*{10.1cm$\Box$ } 
               
                  \noindent \textit{Proof of Theorem $\ref{4.1.}$}.  
             From assumption $\hyperref[f1]{(f_1)}$  and Proposition  $\ref{pro3}$, we have 
             $$
             \begin{aligned}
                    I(u)&=\int_{\Omega}\int_{\Omega}A\left( \dfrac{|u(x)-u(y)|}{|x-y|^s }\right)  \dfrac{dxdy}{|x-y|^{N} } -\int_{\Omega}F(x,u)dx\\
                    &\geqslant\int_{\Omega}\int_{\Omega}A\left( \dfrac{|u(x)-u(y)|}{|x-y|^s }\right)  \dfrac{dxdy}{|x-y|^{N} }-\theta_1 \int_{\Omega} |u|^qdx-\theta_1|\Omega|\\
                    &\geqslant ||u||^{p_0} -\theta_1 C ||u||^q-\theta_1|\Omega|,
                   \end{aligned}
                   $$
          since $p_0>q$, so we have $I(u)\longrightarrow \infty$ as $||u||\longrightarrow \infty$, by Lemma $\ref{4.3}$ $I$ is weakly lower semi-continuous on $W^s_0L_A(\Omega)$, then by Theorem \ref{AA2.2} functional $I$ has a minimum point $u_0$ in $W^s_0L_A(\Omega)$ and $u_0$ is a weakly solution of problem $\hyperref[P_a]{(P_a)}$.
          
          Next we need to verify that $u_0$ is nontrivial. Let $x_0\in \Omega_0$, $0 < R < 1$ satisfy $B_{2R}(x_0)\subset \Omega_0$, where $B_{2R}(x_0)$
          is the ball of radius $2R$ with center at the point $x_0$ in $\R^N$. Let $\varphi\in C_0^\infty(B_{2R}(x_0))$ satisfies $0\leqslant \varphi \leqslant 1$ and
          $\varphi \equiv 1$  in $B_{2R}(x_0)$. Theorem $\ref{3.1..}$ implies that $||\varphi||<\infty.$ Then for $0 < t < 1$, by  $\hyperref[f2]{(f_2)}$, we have    
           $$
             \begin{aligned}
                    I(t\varphi)&= \int_{\Omega}\int_{\Omega}A\left( \dfrac{|t\varphi(x)-t\varphi(y)|}{|x-y|^s }\right) \dfrac{dxdy}{|x-y|^N}  -\int_{\Omega}F(x,t\varphi)dx\\
                    &\leqslant  ||t\varphi||^{\beta}-\int_{\Omega_0}\dfrac{\theta_2}{q} |t\varphi|^{q}dx\\
                    &\leqslant ||\varphi||^{\beta}t^{\beta}-\dfrac{\theta_2}{q}t^q\int_{\Omega_0} |\varphi|^{q}dx\\
                    &\leqslant ct^{\beta}-\dfrac{\theta_2}{q}t^q\int_{\Omega_0} |\varphi|^{q}dx,
                   \end{aligned}
                   $$
          where $\beta=\left\lbrace  p^0 \text{ or } p_0\right\rbrace $ and $c$ is a positive constant. Since $\beta > q$ and $\displaystyle\int_{\Omega_0} |\varphi|^{q}dx>0,$ we have $I(t_0\varphi)<0$ for
          $t_0\in(0,t)$ sufficiently small. Hence, the critical point $u_0$ of functional $I$ satisfies $I(u_0)\leqslant I(t_0\varphi)< 0=I(0)$,
          that is $u_0\neq 0$. 
                  \hspace*{14.2cm $\Box$ }\\          
           \noindent \textit{Proof  of Corollary $\ref{cr}$}. In view of the proof of Theorem $\ref{4.1.}$, we only need to check that $I(u)\longrightarrow \infty$ as $||u||\rightarrow \infty$. Let $u\in W^s_0L_A(\Omega)$ with $||u||>1$, since
           $p_0 = q$ and $\theta_1< (\lambda_1)/(2)$, by assumptions $\hyperref[f1]{(f_1)}$ and  ($\ref{lam}$), we have
               $$
              \begin{aligned}
                     I(u)&=\int_{\Omega}\int_{\Omega}A\left( \dfrac{|u(x)-u(y)|}{|x-y|^s }\right)  \dfrac{dxdy}{|x-y|^{N} } -\int_{\Omega}F(x,u)dx\\
                     &\geqslant \int_{\Omega}\int_{\Omega}A\left( \dfrac{|u(x)-u(y)|}{|x-y|^s }\right)  \dfrac{dxdy}{|x-y|^{N} }-\theta_1 \int_{\Omega} |u|^{p_0}dx-\theta_1|\Omega|\\
                     &\geqslant ||u||^{p_0} -\theta_1 \frac{1}{\lambda_1} ||u||^{p_0}-\theta_1|\Omega|,\\
                     & =\left( 1 -\theta_1\frac{1}{\lambda_1}\right) ||u||^{p_0}-\theta_1|\Omega|.
                    \end{aligned}
                    $$
                    So we have $I(u)\longrightarrow \infty$ as $||u||\longrightarrow \infty$.
                         \hspace*{7.2cm$\Box$ } 
\begin{remark}
Evidently, if $f(x, 0) \neq 0$ a.e. in $\Omega$, then the weak solution obtained in Corollary $\ref{cr}$ is nontrivial.
\end{remark}         
\bibliographystyle{amsplain}

\end{document}